\documentclass[12pt]{article} 
\def\Bbb{\bf} 
\newcommand\C{{ \Bbb C}}
\newcommand\R{{\Bbb R}}
\newcommand\Z{{\Bbb Z}}
\newcommand\PP{{\Bbb P}}

\newcommand\G{{\Gamma}} 

\newcommand\e{{\varepsilon}}
\newtheorem{lmm}{Lemma}

\newtheorem{rmk}{Remark}
\newtheorem{prp}{Proposition}

\def\comment#1{ }

\newcommand{\dfrac}[2]{%
\frac{\displaystyle{#1}}{\displaystyle{#2}}}
\begin{document}
\title{
A family of Schottky groups arising from the hypergeometric equation (tex/schot/IY2)}
\author{Takashi ICHIKAWA \thanks{Department of Mathematics, Faculty of Science and Engineering, Saga University, Saga 840-8502 Japan, 
{\it E-mail address}: ichikawa@ms.saga-u.ac.jp} 
\and Masaaki YOSHIDA \thanks{Department of Mathematics, Kyushu University, 
                             Fukuoka 810-8560 Japan, 
{\it E-mail address}: myoshida@math.kyushu-u.ac.jp}}
\maketitle

\noindent
{\textbf{Abstract:}}\quad
We study a complex 3-dimensional family of classical Schottky groups of 
genus 2 as monodromy groups of the hypergeometric equation. 
We find non-trivial loops in the deformation space; these correspond to 
continuous integer-shifts of the parameters of the equation.
{\textbf{Running title:}}\quad Schottky groups arising from the hypergeometric equation\\
\tableofcontents
\section{Introduction}
When the three exponents of the hypergeometric differential equation are 
purely imaginary, its monodromy group is a classical Schottky group of 
genus 2; such groups form a real 3-dimensional family. 
Since, under a small deformation, a Schottky group remains to be a 
Schottky group, by deforming the parameters of the hypergeometric equation, 
we have a complex 3-dimensional family of Schottky groups. 
We introduce (in \S4) a complex 3-dimensional family $S$ of classical Schottky groups, containing the above ones coming from pure-imaginary-exponts cases, equipped with some additional structure. We study a structure of $S$, and construct loops (real 1-dimensional families of classical Schottky groups) in $S$ generating the fundamental group of $S$. These loops correspond to {\it continous} interger-shifts of the parameters of the hypergeometric equation.
\section{The hypergeometric equation}
We consider the hypergeometric differential equation
$$E(a,b,c): x(1-x)\frac{d^2u}{dx^2}+\{c-(a+b+1)x\}\frac{du}{dx}-abu=0.$$
For (any) two linearly independent solutions $u_1$ and $u_2$, the (multi-valued) map
$$s:X:=\C-\{0,1\}\ni x\longmapsto z=u_1(x):u_2(x)\in \PP^1:=\C\cup\{\infty\}$$
is called a Schwarz map (or Schwarz's $s$-map). 
If we choose as solutions $u_1$ and $u_2$, near $x=0$,  a holomorphic one  and
$x^{1-c}$ times a holomorphic one, respectively, then the circuit matrices
around $x=0$ and  $1$ are given by
$$\gamma_1=\left(\begin{array}{cc}1&0\\0&e^{2\pi i(1-c)}\end{array}\right)
\quad \mbox{and}\quad
\gamma_2=P^{-1}\left(\begin{array}{cc}1&0\\0&e^{2\pi i(c-a-b)}\end{array}\right)P,$$
respectively, where $P$ is a connection matrix given by
$$P=\pmatrix{
\displaystyle{\G(c)\G(c-a-b)\over\G(c-a)\G(c-b)}&
\displaystyle{\G(2-c)\G(c-a-b)\over\G(1-a)\G(1-b)}\cr
\displaystyle{\G(c)\G(a+b-c)\over\G(a)\G(b)}&
\displaystyle{\G(2-c)\G(a+b-c)\over\G(a-c+1)\G(b-c+1)}};
$$
here $\G$ denotes the Gamma function.
Note that the matrices act the row vector $(u_1,u_2)$ from the right.
These generate the monodromy group $M(a,b,c)$ of the equation $E(a,b,c)$.
\section{Hypergeometric equations with purely imaginary exponent-differences ([SY], [IY])}
If the exponent-differences
$$\lambda=1-c,\quad \mu=c-a-b,\quad \nu=b-a$$
(at the singular points $x=0,1$ and $\infty$, respectively) are purely imaginary, the image of the upper-half part
$$X^+:=\{x\in X\mid \Im x\ge0\}$$
is bounded by the three circles (in the $s$-plane) which are images of the three intervals $(-\infty,0),(0,1)$ and $(1,+\infty)$; put
$$C_1=s((-\infty,0)),\quad C_3=s((0,1)),\quad C_2=s((0,+\infty)).$$
Note that if we continue analytically $s$ through the interval, say $(0,1)$, to the lower-half part $X^-:=\{x\in X\mid \Im x\le0\}$, then the image of $X^-$ is the mirror image of that of $X^+$ under the reflection with center $C_3$.
\par\medskip\noindent
The reflection with center (mirror) $C_j$ will be denoted by $\rho_j\ (j=1,2,3)$. The monodromy group $M(a,b,c)$ is the group consisting of even words of $\rho_j\ (j=1,2,3)$. This group is a Schottky group of genus 2. The domain of discontinuity modulo $M(a,b,c)$ is a curve of genus 2 defined over reals. 
\par\medskip\noindent
Choosing the solutions $u_1$ and $u_2$ suitably, we can assume that the
centers of the three circles are on the real axis, and that the circles
$C_2$ and $C_3$ are inside the circle $C_1$ (see Figure 1). 
\begin{figure}[hhh]
\begin{center}\include{schot0}\end{center}
\caption{The circles and the disks}
\end{figure}
Let $C_1'$
and $C_2'$ be the mirror images with respect to $C_3$ of the circles
$C_1$ and $C_2$, respectively. Then $\gamma_1:=\rho_3\rho_1$ maps inside
of $C_1$ onto inside of $C_1'$; and $\gamma_2:=\rho_3\rho_2$ maps
inside of $C_2$ onto the outside of $C_2'$. Note that one of the two fixed points of $\gamma_1$ is inside $C_1'$,  and the other one is outside $C_1$; one of the two fixed  points of $\gamma_2$ is inside $C_2'$, and the other one is inside $C_2$. The transformations $\gamma_1$ and $\gamma_2$ are conjugate to those
given in \S2.
\par\medskip
We thus have four disjoint disks
$D_1,D_1',D_2,$ and $D_2',$ whose centers are on the real axis, and two
loxodromic transformations $\gamma_1$ and $\gamma_2$ taking $D_1$ and
$D_2$ to the complementary disks of $D_1'$ and $D_2'$, respectively. Note that one of the two fixed points of $\gamma_1$ is in $D_1'$ and the other one in $D_1$; one of the two fixed points of $\gamma_2$ is in $D_2'$ and the other one in $D_2$.
\section{The moduli space $S$}
Two loxodromic transformations $\gamma_1$ and
$\gamma_2$ generate a (classical) Schottky group if and only if there
are four disjoint closed disks $D_1,D_1',D_2,$ and $D_2',$ such that 
$\gamma_1$ and $\gamma_2$ take $D_1$ and
$D_2$ to the complementary disks of $D_1'$ and $D_2'$, respectively. Throughout the
paper, disks are always assumed to be {\it closed};
we simply call them {\it disks}. The closure of the complement of a disk in $\PP^1$  will simply be called {\it the complementary disk of the disk}.
\par\smallskip\noindent
{\bf Definition}\ \ {\it
Let $S$ be the space of two loxodromic transformations $\gamma_1$ and
$\gamma_2$ equipped with four disjoint disks $D_1,D_1',D_2,$ and $D_2',$ 
such that $\gamma_1$ and $\gamma_2$ take $D_1$ and
$D_2$ to the complementary disks of $D_1'$ and $D_2'$, respectively.} 
\par\smallskip\noindent
The space has a natural structure of complex 3-dimensional manifold.
We are interested in its homotopic property. 
\par\medskip
Note that if a loxodromic transformation $\gamma$ takes a disk $D$ onto 
the complementary disk of a disk $D'$ ($D\cap D'=\emptyset$), then $\gamma$ has 
a fixed point in $D$ and the other fixed point in $D'$.
\par\smallskip\noindent
A loxodromic transformation is determined by the two fixed points and
the multiplier.
\par\smallskip\noindent
\begin{lmm} 
For given two disjoint disks $D$ and $D'$, there is a unique point
$F\in D$ (resp. $F'\in D'$) such that by a(ny) linear fractional 
transformation taking $F$ (resp. $F'$) to $\infty$, the two circles 
$\partial D$ and $\partial D'$ are transformed into cocentric circles.
\end{lmm}
\par\smallskip\noindent
Proof. Let the two disks be given as
$$D:|z-a|\le r,\quad D':|z-a'|\le r'.$$
By the transformation $z\to w$ defined by
$$z=\frac1w+\zeta$$
taking $\zeta$ to $\infty$, the circles $C=\partial D$ and $C'=\partial D'$
are transformed into circles with centers
$$c:=\frac{\bar \zeta-\bar a}{r^2-|\zeta-a|^2},\quad 
 c':=\frac{\bar \zeta-\bar a'}{r'^2-|\zeta-a'|^2},$$
respectively. Equating $c$ and $c'$, we have the quadric equation
$$\zeta^2+\left(-a-a'+\frac{r^2-r'^2}{\bar a-\bar{a'}}\right)\zeta+aa'+\frac{ar'^2-a'r^2}{\bar a-\bar{a'}}=0.$$
Put $\zeta=(a-a')\eta+a'$. Note that $\zeta=a'$ and $a$ correspond to
$\eta=0$ and $1$, respectively. Then the equation above for $\zeta$ reduces to
$$\eta^2+\left(-1+\frac{r^2-r'^2}{|a-a'|^2}\right)\eta+\frac{r'^2}{|a-a'|^2}=0.$$ It is a high school mathematics to see that each of the two roots of this equation is in each of the two intervals
$$\left(0,\frac{r'}{|a-a'|}\right)\mbox{\quad and\quad} \left(1-\frac{r}{|a-a'|},1\right).$$
Since cocentric circles are mapped to cocentric circles under linear transformations, this completes the proof.\qquad []
\begin{rmk} The circles $\partial D$ and $\partial D'$ are Apollonius circles with respect to the two centers $F$ and $F'$.\end{rmk}
\begin{lmm} 
For given two disjoint disks $D$ and $D'$ and a point $f'$ in 
the interior of $ D'$, there is a 1-parameter family 
(parametrized by a circle) of loxodromic 
transformations $\gamma$ which take $D$
onto the complementary disk of $D'$, and fix $f'$. 
Absolute value $|m|$ of the multiplier 
$m$ of $\gamma$ is determined by the given data. 
\par\noindent
{\rm (1)} If $f'\not=F'$, then the other fixed point $f\in D$ of $\gamma$ must be on 
the Apollonius circle $A=A(f')$ in $D$ determined by the two centers of the circles $\partial D$ and $\partial D'$, and propotion $|m|$. 
$\arg m\in\R/2\pi\Z$ determines $f$, and vice versa. 
\par\noindent
{\rm (2)} If $f'=F'$, then the other fixed point is $F\in D$. $\arg m\in \R/2\pi\Z$ remains free.
\end{lmm}
\par\smallskip\noindent
Proof. We can assume that $f=0$ and $f'=\infty$, so that the
transformation in question can be presented by $z\mapsto mz.$
Let $c$ and $r$ be the center and the radius of the disk $D$, and $c'$ and 
$r'$ those of the complementary disk of $D'$. Then we have $c'=mc$ and
$r'=|m|r$. 
\par\noindent
(1) If $c\not=c'$, $f(=0)$ is on the Apollonius circle
$$A:|f-c'|=|m||f-c|$$
with centers $c$ and $c'$, and proportion $m$.  
It is easy to see that this circle is in $D$. (See Figure 2.)
\begin{figure}[hhh]
\begin{center}\include{apollo}\end{center}
\caption{Apollonius circle in $D$}
\end{figure}
\par\noindent
(2) If $c=c'=0$, then the Apollonius circle reduces to a point.\qquad []
\par\medskip\noindent
In the case (2), we blow up the point $F$ to be the circle $\R/2\pi\Z$, and call this circle also the Apollonius circle $A$. Under this convention, two disjoint disks $D$ and $D'$, an interior point $f'\in D'$, and a point $f$ on the Apollonius circle $A$ determines uniquely a loxodromic transformation.
\par\smallskip\noindent
Thus an element of $S$ can be determined by four disjoint disks 
$D_1,D_1',D_2,D_2',$ two interior points $f_1'\in D_1', f_2'\in D_2'$, 
and two points $f_1$ on the Apollonius circle $A_1$ (determined by 
$D_1,D_1', f_1'$) in $D_1$,  
and $f_2$ on the Apollonius circle $A_2$ (determined by  $D_2,D_2', f_2'$)
in $D_2$. 
Since disks are contractible, we have
\begin{prp} The fundamental group of $S$ can be generated by 
the moves of the four disjoint disks $D_1,D_1',D_2,D_2',$ 
and the moves of $f_1$ in $A_1$, and $f_2$ in $A_2$.
\end{prp}
In the next section we explicitly construct loops (real 1-dimensional families of classical Schottky groups) in $S$.
\section{Loops in $S$}
Let $\gamma_1$ and $\gamma_2$ be as in \S2. The fixed points of
$\gamma_1$ are $\{0,\infty\}$, and those of $\gamma_2$ are
$$f_2=\frac{\G(c)\G(a-c+1)\G(b-c+1)}{\G(2-c)\G(a)\G(b)}\quad\mbox{and}\quad
f_2'=\frac{\G(c)\G(1-a)\G(1-b)}{\G(2-c)\G(c-a)\G(c-b)}.$$
We change the coordinate $z$ by multiplying $1/f_2$. Then the fixed points of
$\gamma_1$ remain to be  $\{0,\infty\}$, 
and those of $\gamma_2 $ become $1$ and
$$\alpha=g(a)g(b),\quad\mbox{where}\quad 
g(x)=\frac{\sin(\pi c-\pi x)}{\sin(\pi x)}.$$
Indeed we have
$$\frac{f_2'}{f_2}=\frac{\G(1-a)\G(1-b)}{\G(c-a)\G(c-b)}\cdot
\frac{\G(a)\G(b)}{\G(a-c+1)\G(b-c+1)}\quad\mbox{and}\quad
\G(x)\G(1-x)=\frac\pi{\sin(\pi x)}.$$
\par\medskip\noindent
{\bf Definition}\ \ {\it When the three exponent-differences are purely
imaginary:
$$\lambda=i\theta_0,\quad \mu=i\theta_1,\quad \nu=i\theta_2,\quad \theta_0,\theta_1,\theta_2>0,$$
the generators $\{\gamma_1,\gamma_2\}$ of the monodromy group of
$E(a,b,c)$ given in \S2 (take the four disks given in \S3) form a simply
connected real 3-dimensional submanifold of $S$; this submanifold is called $S_0$.}
\par\medskip\noindent
In this section, we construct loops in $S$, with base in $S_0$, which
generate the fundamental group of $S$. When the exponent-differences are
as above, note that the parameters can be expressed as
$$a=\frac12-\frac{i}2(\theta_0+\theta_1+\theta_2),\quad 
b=\frac12-\frac{i}2(\theta_0+\theta_1-\theta_2),\quad c=1-i\theta_0.$$
For clarity, we denote the four disks $D_1,D_1',D_2,D_2'$ by
$$D_0\ (\ni0),\quad D_\infty\ (\ni\infty),\quad D_1\ (\ni 1),\quad
D_\alpha\ (\ni\alpha).$$
\subsection{The disk $D_\alpha$ travels around the disk $D_0$}
We construct a loop in $S$ (with base in $S_0$), which is represented by 
a travel of the disk $D_\alpha$  around the disk $D_0$ with a change of argument (by $2\pi$) of the multiplier of $\gamma_2$. 
We fix $c$ and the real part of $a$ (as $1/2$), and let the real part of 
$b$ move from $1/2$ to $3/2$; the imaginary parts of $a$ and $b$ are so
chosen that the monodromy group $M(a,b,c)$ remains to be a Schottky
group along the move. Putting $c=1-i\theta_0$, we have
$$g(x)=\dfrac{\e^{-1}-\e e^{2\pi ix}}{1-e^{2\pi ix}}=\e+\frac{\e^{-1}-\e}{1-e^{2\pi ix}}, \quad\mbox{where}\quad \e=e^{-\pi\theta_0}<1.$$
Set
$$r=g\left(\frac12-\frac{i}2(\theta_0+2)\right),\quad 
R=g\left(\frac12-\frac{i}2(\theta_0+1)\right).$$
Then we have 
$$0<\e<r=\e\frac{e^{2\pi}+\e^{-1}}{e^{2\pi}+\e}<
R=\e\frac{e^{\pi}+\e^{-1}}{e^{\pi}+\e}<1.$$
Let $D_0$ be the disk with center at $0$ and with radius $\e r$, and $D_\infty$ the complementary disk of the disk with center at $0$ and with radius $\e^{-1}r$. Note that we have
$$0<\e r<\e R<R<1<\e^{-1}r,$$
and that $\gamma_1$ maps $D_0$ onto the outside disk of $D_\infty$.

Define real continuous functions $\phi(t)$ and $\psi(t)$ for $0\le t\le1$ by
$$Re^{2\pi it}=g\left(\frac12+\phi(t)-\frac{i}2(\theta_0+\psi(t))\right).$$
Since 
$$2\pi i\left(\frac12+\phi(t)-\frac{i}2(\theta_0+\psi(t))\right)=\log(Re^{2\pi it}-\e^{-1})-\log(Re^{2\pi it}-\e),$$
the function $\phi(t)$ is monotone increasing with $\phi(0)=0, \phi(1)=1$,
and $\psi(t)$ satisfies $\psi(0)=\psi(1)=1.$ Choose $\theta_1$ so big that
\begin{itemize}
\item $\theta_1>\max\{\psi(t)\mid0\le t\le1\}$, and that
\item for any $z$ in the ring $\{\e R\le|z|\le R\}$, there are two disjoint disks $D_1\ (\ni1)$ and $D_z\ (\ni z)$, and a fractional linear transformation with multiplier $e^{-2\pi\theta_1}$ mapping $D_1$ onto the complementary disk of $D_z$. 
\end{itemize} 
\begin{figure}[hhh]
\begin{center}\include{travel3}\end{center}
\caption{The disk $D_\alpha$ travels around the disk $D_0$}
\end{figure}
Now set $\theta_2(t)=\theta_1-\psi(t)>0$ and deform the parameters $a$ and $b$ as
$$\begin{array}{rl}
a(t)&=\dfrac12-\dfrac{i}2(\theta_0+\theta_1+\theta_2(t)),\\[2mm]
b(t)&=\dfrac12+\phi(t)-\dfrac{i}2(\theta_0+\theta_1-\theta_2(t)).
\end{array}$$
Then we have $\e<g(a(t))<1, |g(b(t))|=R$, 
and so $\alpha(t)=g(a(t))g(b(t))$ satisfies
$$\e R<|\alpha(t)|=g(a(t))|g(b(t))|<R.$$
Thus $\alpha(t)$ travels around the disk $D_0$ in the ring 
$\{\e R\le|z|\le R\}$, 
and there are two disjoint disks $D_1\ (\ni1)$ and $D_{\alpha(t)}\ (\ni
\alpha)$ in the ring $\{\e r<|z|<\e^{-1}r\}$, and a transformation $\gamma_2(t)$ with multiplier $e^{2\pi i\mu}=e^{-2\pi\theta_1}$ which maps $D_1$ onto the complementary disk of $D_{\alpha(t)}$. 
\subsection{The disk $D_\alpha$ travels around the disk $D_1$}
We construct a loop in $S$ (with base in $S_0$), which is represented by 
a travel of the disk $D_\alpha$  around the disk $D_1$ with a change of argument (by $2\pi$) of the multiplier of $\gamma_2$. 
We fix $c,a$ and the imaginary part of $b$, and let the real part of 
$b$ move from $1/2$ to $3/2$. Take $\theta_0>0$ and $\theta'<0$ satisfying

\begin{equation}
\frac{e^{\pi \theta'} (1 + e^{\pi \theta'})}{1 - e^{\pi \theta'}} 
< \varepsilon^{-2}, \quad \varepsilon = e^{- \pi \theta_0} . 
\end{equation}
Set 
$$
r = \varepsilon + \frac{\varepsilon^{-1} - \varepsilon}{1 + e^{\pi \theta'}}, 
\quad 
R = \varepsilon + \frac{\varepsilon^{-1} - \varepsilon}{1 - e^{\pi \theta'}}. 
$$
Since we have 
$$ 
\frac{1 - \varepsilon r }{r} = \frac{1 - \varepsilon^2}{r} 
\left( 1 - \frac{1}{1 + e^{\pi \theta'}} \right) > 0 
$$
and 
$$
\varepsilon^{-1} r - \varepsilon R = 
(1 - \varepsilon^2) \left( 1 + \frac{\varepsilon^{-2}}{1 + e^{\pi \theta'}} - 
\frac{1}{1 - e^{\pi \theta'}} \right) > 0, 
$$
where the last inequality holds thanks to (1), 
there is a positive number $s$ satisfying 
$$
0 < s < \min \left\{ r, \ \frac{1 - \varepsilon r }{r}, \ 
\frac{\varepsilon^{-1} r - \varepsilon R }{\varepsilon^{-1} + R} \right\}, 
$$
that implies 
\begin{equation}
(\varepsilon + s) r < 1, \quad (\varepsilon + s) R < \varepsilon^{-1} (r - s). 
\end{equation}
Take $\theta_1$ so big that 
\begin{equation}
\varepsilon + 
\frac{\varepsilon^{-1} - \varepsilon}{1 + e^{\pi(\theta_0 + \theta_1)}} 
< \varepsilon + s. 
\end{equation}
Now set $\theta_2 := \theta_0 + \theta_1 - \theta'$, and 
$$
a = \frac{1}{2} - \frac{i}{2} ( \theta_0 + \theta_1 + \theta_2 ), 
$$
and deform (the real part of) $b$ as
$$
b(t) = \frac{1}{2} + t - \frac{i}{2} ( \theta_0 + \theta_1 - \theta_2 ) 
= \frac{1}{2} + t - \frac{i}{2} \theta' . 
$$
Note that 
$$
\varepsilon < g(a)  =  \varepsilon + \frac{\varepsilon^{-1} - \varepsilon}
{1 + e^{\pi(\theta_0 + \theta_1 + \theta_2)}} < \varepsilon + s, 
\quad \mbox{here we used (3),}$$ and that  
$$           g(b(t))  =  \varepsilon + \frac{\varepsilon^{-1} - \varepsilon}
{1 + e^{\pi \theta'} e^{2 \pi i t}} \qquad 
\mbox{satisfies} \qquad r \leq |g(b(t))| \leq R.$$
These inequalities together with (2) implies that 
$\alpha(t) = g(a) g(b(t))$ is in the ring
$\{ \varepsilon r \leq |z| \leq (\varepsilon + s) R \}$ and that
$$
\alpha(0) = \alpha(1) \leq (\varepsilon + s) r < 1 < \varepsilon R 
\leq \alpha(1/2).$$
This shows that $\alpha(t)$ travels around $1$. 
\begin{figure}[hhh]
\begin{center}\include{travel4}\end{center}
\caption{The disk $D_\alpha$ travels around the disk $D_1$}
\end{figure}
Let $D_0$ be the disk with center at $0$ and with radius $\varepsilon (r-s)$, 
$D_{\infty}$ the complementary disk of the disk with center at $0$ and 
with radius $\varepsilon^{-1} (r-s)$. 
Then $\gamma_1$ maps $D_0$ onto the complementary disk of $D_{\infty}$, and by right (2), 
$D_0 \cup D_{\infty}$ is disjoint from the ring 
$$
\{ \varepsilon r \leq |z| \leq (\varepsilon + s) R \} \ni 1, \alpha(t). 
$$ 
Thus by taking $\theta_1 > 0$ sufficiently big, for any $0 \leq t \leq 1$, 
there are two disjoint disks 
$D_1 \ (\ni 1)$, $D_{\alpha(t)} \ (\ni \alpha(t))$ 
in the complement of $D_0 \cup D_{\infty}$ 
and a transformation $\gamma_2(t)$ with multiplier 
$e^{2 \pi i \mu} = e^{-2 \pi \theta_1} e^{-2 \pi i t}$ 
which maps $D_1$ onto the complementary disk of the $D_{\alpha(t)}$. 

\subsection{The multiplier of $\gamma_2$ travels around $0$}
We construct a loop in $S$ (with base in $S_0$), which is represented by 
the change of argument (by $2\pi$) of the multiplier of $\gamma_2$. 
We fix $b$ and $c$, and the imaginary part of $a$, and let the real part 
of $a$ move from $1/2$ to $3/2$.
Set
$$\theta:=\theta_0+\theta_1+\theta_2,\quad \e:=e^{-\pi\theta_0},\quad 
r:=\frac{\e^{-1}-\e}{e^{\pi\theta}-1}.$$
Choose and fix $\theta_0,\theta_1$ and $\theta_2$ so that 
$\theta_1=\theta_2$ and
$$r<\min\left\{1-\e,\ \frac1{1+\e^{-1}}\right\}.$$
Then since
$$b=\frac12-\frac{i}2(\theta_0+\theta_1-\theta_2)=\frac12-\frac{i}2\theta_0,$$
we have $e^{2\pi ib}=-\e^{-1}$, and so $g(b)=1$. 
\par\smallskip
Now we deform the parameter $a$ as
$$a(t)=\frac12+t-\frac{i}2\theta,\quad 0\le t\le1.$$
Since we have
$$\alpha(t)=g(a(t))g(b)=g(a(t))=\e+\frac{\e^{-1}-\e}{1+e^{\pi\theta}e^{2\pi it}},$$
the point $\alpha$ is in the disk with center at $\e$ and with radius $r$.
We thus have
$$
\varepsilon^2 (1 + r) < \varepsilon - r \leq |\alpha(t)| \leq \varepsilon + r 
< 1 < 1 + r. 
$$
\begin{figure}[hhh]
\begin{center}\include{travel5}\end{center}
\caption{The fixed point $\alpha$ of $\gamma_2$ travels along the Apollonius circle}
\end{figure}
\noindent
Let $D_0$ be the disk with center at $0$ and with radius $\varepsilon^2(1+r)$, 
$D_{\infty}$ the complementary disk of the disk with center at $0$ and with radius $1+r$, 
$D_{\alpha(t)} \ (\ni \alpha(t))$ the disk with center at $\varepsilon$ and 
with radius $r$. 
Then $\gamma_1$ maps $D_0$ onto the complementary disk of $D_{\infty}$, 
and $\alpha(t)$ belongs to the disk $\{ | z - \varepsilon | \leq r \}$ 
which is disjoint from $D_0 \cup D_{\infty} \cup \{ 1 \}$. 
Thus by taking $\theta_1 = \theta_2 > 0$ sufficiently big, 
for any $0 \leq t \leq 1$, there are two disjoint disks 
$D_1 \ (\ni 1)$, $D_{\alpha(t)} \ (\ni \alpha(t))$ 
in the complement of $D_0 \cup D_{\infty}$ 
and a transformation $\gamma_2(t)$ with multiplier 
$e^{2 \pi i \mu} = e^{-2 \pi \theta_1} e^{-2 \pi i t}$ 
which maps $D_1$ onto the complementary disk of $D_{\alpha(t)}$. 

\subsection{The multiplier of $\gamma_1$ travels around $0$}
We construct a loop in $S$ (with base in $S_0$), which is represented by 
the change of argument (by $2\pi$) of the multiplier of $\gamma_1$.
Note that by the change of variable $x\to1-x$, the equation $E(a,b,c)$ changes into $E(a,b,a+b+1-c)$. So we have only to literally follow \S5.3 exchanging $c$ and $a+b+1-c$, and $\theta_0$ and $\theta_1$. We thus fix $b$ and the imaginary parts of $a$ and $c$, and let the real part of $c$ move from 1 to 2, and let the real part of $a$ move from $1/2$ to $3/2$ keeping $a+b+1-c$ constant.
\section{Miscellanea}
For a Schottky group $\Gamma$ of genus 2, the quotient of the domain of discontinuity modulo $\Gamma$ is a curve of genus 2. A curve of genus 2 is a double cover of $\PP^1$ branching at six points, which are uniquely determined modulo automorphisms of $\PP^1$ by the curve. Thus a Schottky group determines a point of the configuration space $X\{6\}$ of six-point sets on $\PP^1$. (The fundamental group of $X\{6\}$ is the Braid group with five strings.) When  all the exponent-differences of the hypergeometric equation are real, its monodromy group is a Schottky group, which determines six points on a line. Thus they determines a point of the configuration space $X(6)$ of colored six point on $\PP^1$. (The fundamental group of $X(6)$ is the colored Braid group with five strings.) The space $X(6)$ is well-studied.  
\begin{rmk}When all the exponent-differences are real, there is a Schottky automorphic function defined by an absolutely convergent infinite product, which induces a holomorphic map of the genus 2 curve onto $\PP^1$ {\rm([IY])}. This infinite product remains to be convergent for groups represented by loops above, if we take the multipliers of $\gamma_1$ and $\gamma_2$ sufficiently small. This is because there is a circle separating two disks among the four {\rm([BBEIM, Chapter 5])}. \end{rmk}
\par\smallskip\noindent
{\bf Problems:}\quad The four loops constructed in \S5 induces those in $X\{6\}$. Do they generate the fundamental group of $X\{6\}$? Same problem for $X(6)$. 
\par\bigskip\noindent
{\textbf{Acknowledgement:}}\quad
The authors are grateful to Professors Susumu Hirose, Ken'ichi Ohshika, Hiroki Sato and Hiro-o Yamamoto for their encouragement.


\begin{thebibliography}{BBEIM}
\bibitem[BBEIM]{BBEIM} 
{\sc E.\,D.\,Belokolos, A.\,I.\,Bobenko, V.\,Z.\,Enol'skii, A.\,R.\,Its and V.\,B.\,Matveev, }
{\it Algebro-geometric Approach to Nonlinear Integrable Equations,} 
Springer Series in Nonlinear Dynamics 
(Springer-Verlag, 1994). 
%
%
\bibitem[IY]{IY}{\sc T. Ichikawa and M.\,Yoshida,} 
On Schottky groups arising from the hypergeometric equation with imaginary exponents, Proc AMS 132(2003), 447--454.
%
%
\bibitem[SY]{SY}{\sc T.\,Sasaki and M.\,Yoshida,} 
A geometric study of the hypergeometric fucntion with imaginary exponents, Experimental Math. {\bf 10}(2000), 321 --330.
%
\bibitem[Sch]{Sch} 
{\sc F.\,Schottky, }
\"{U}ber eine specielle Function, welche bei einer bestimmten 
linearen Transformation ihres Arguments univer\"{a}ndert bleibt, 
J. Reine Angew. Math. 101 (1887) 227-272. 
\end{thebibliography}
\end{document}